%% file: 0paper.tex
\documentclass[12pt]{article}
\usepackage{amsfonts}
\usepackage{epsfig}
\title{Lenghtening a Tetrahedron}
\author{Richard Evan Schwartz \thanks{\hskip 5 pt Supported by 
N.S.F. Research Grant DMS-1204471}}

\newtheorem{theorem}{Theorem}[section]

\newtheorem{lemma}[theorem]{Lemma}

\newtheorem{corollary}[theorem]{Corollary}

\def\startproof{{\bf {\medskip}{\noindent}Proof: }}

\def\endproof{$\spadesuit$  \newline}

\def\A{\mbox{\boldmath{$A$}}}%
\def\B{\mbox{\boldmath{$B$}}}%
\def\C{\mbox{\boldmath{$C$}}}%
\def\D{\mbox{\boldmath{$D$}}}%
\def\N{\mbox{\boldmath{$N$}}}%
\def\P{\mbox{\boldmath{$P$}}}%
\def\R{\mbox{\boldmath{$R$}}}%
\def\Z{\mbox{\boldmath{$Z$}}}%

\begin{document}
\maketitle
\begin{abstract}
We give rigorous, computer assisted
proofs of a number of statements about
the effect on the volume 
of lengthening various edges of
a tetrahedron.  Our results give new
and sharp polynomial inequalities concerning
the Cayley-Menger determinant and its partial
derivatives.
\end{abstract}

\input{1intro}
\input{2setup}
\input{3posdom}
\input{4code}
\input{5select}
\input{6exist}

\input{refs}

\end{document}

%% file: 1intro.tex
\section{Introduction}

This paper was inspired by a question posed by Daryl Cooper:
{\it Suppose that you lengthen all the sides of a 
tetrahedron by one unit. Is the result still a tetrahedron, 
and (if so) does the volume increase?\/}
More formally, say that a list
$\{d_{ij}|\ i \not = j \in \{1,2,3,4\}\}$ 
is {\it tetrahedral\/} if there are $4$ distinct points
$V_1,V_2,V_3,V_4 \in \R^3$ so that
$d_{ij}=\|V_i-V_j\|$ for all $i,j$.
We call the list $\{d_{ij}+1\}$ the {\it unit lengthening\/}
of $\{d_{ij}\}$.

\begin{theorem}
\label{main}
The unit lengthening of a tetrahedral list is also tetrahedral.
If $\Delta_0$ is the original tetrahedron and $\Delta_1$ is
the new tetrahedron, then
$$\frac{{\rm volume\/}(\Delta_1)}{{\rm volume\/}(\Delta_0)} 
\geq \bigg(1+\frac{6}{\sum_{i<j} d_{ij}}\bigg)^3.$$
\end{theorem}
The inequality is sharp, because it is an equality for all
regular tetrahedra.

We also have the following general result.
\begin{theorem}
\label{drs}
\label{exist}  In every dimension,
the unit lengtening of a simplicial list
is again simplicial, and the new simplex
has volume larger than the original.
\end{theorem}
Here a {\it simplicial\/} list is the obvious
generalization of a tetrahedral list to higher dimensions.
One could say that Theorem \ref{drs} is new, and one
could say that it has been there all along.
After discussing an earlier version of this
paper with Peter Doyle and Igor Rivin, they realized
that the general result follows from a theorem,
[{\bf WW\/}, Corollary 4.8], attributed to Von Neumann.
I'll give the argument in an appendix.  It is
independent from the rest of the paper.

Theorem \ref{main} relies on a sharp
inequality concerning the
Cayley-Menger determinant and one of
its directional derivatives.  
Let $K_4$ be the complete graph on $4$ vertices.
Say that a {\it pseudo-tetrahedron\/} 
 is a non-negative  labeling
of the edges of $K_4$ so that, going around any
$3$-cycle of $K_4$, the edges satisfy the triangle
inequality. Let $X$ denote the space of
pseudo-tetrahedra.  We think of $X$ as a polyhedral
cone in $\R^6$ by considering the points
$(d_{12},d_{13},d_{14},d_{23},d_{24},d_{34})$.

Given a pseudo-tetrahedron $D=\{d_{ij}\}$ we have the famous
{\it Cayley-Menger determinant\/}
\begin{equation}
\label{cm}
f(D)=\det\left[\matrix{0&1&1&1&1 \cr
1&0&d_{12}^2& d_{13}^2& d_{14}^2\cr
1&d_{21}^2&0&d_{23}^2& d_{24}^2 \cr
1&d_{31}^2& d_{32}^2&0& d_{34}^2\cr
1&d_{41}^2& d_{42}^2& d_{43}^2&0}\right]
\end{equation}
When $D$ represents a tetrahedron $T_D$, we have
the following classic result.
\begin{equation}
\label{classic}
f(D)=288V^2=2^3 \times (3! V)^2,\hskip 30 pt
V={\rm volume\/}(T_D).
\end{equation}
See [{\P\/}] for a proof, and [{\bf Sa\/}] for a vast survey of
generalizations.
We also define the directional derivative
\begin{equation}
\label{df}
g=D_{(1,1,1,1,1,1)}f.
\end{equation}
$f$ is a homogeneous polynomial of degree $6$ and
$g$ is a homogeneous polynomial of degree $5$.
Theorem \ref{main} is a quick consequence of the
following result.

\begin{theorem}
\label{ineq1}
Let $C$ be a constant.
The function
$g \sum_{i<j} d_{ij} - C f$
is non-negative on $X$ if and only if
$C \in [16,36]$.
\end{theorem}

We will reduce Theorem \ref{ineq1} to the statement that
a certain polynomial in
$\Z[X_1,...,X_5]$ is non-negative on the unit 
cube $[0,1]^5$.  We then use about
an hour of exact integer calculation in Java to establish
the non-negativity.  \S 3 describes the method and
\S 4 gives details about its implementation in this case. I call it
the Method of Positive Dominance.  
I have no idea if it is a known technique,
though I also used it in [{\bf S\/}].
\newline

In \S 5, we will use the same methods to prove 
generalizations of Theorems \ref{main} and \ref{ineq1}
which deal with selectively lengthening some subset
of the edges of a tetrahedron.
Here is the framework for these results.  Each
pseudo-tetrahedron gives rise to $4$ {\it vertex sums\/}
and $3$ {\it axis sums\/}.  A vertex sum is the sum
of the labels of $3$ edges incident to a given vertex -- e.g.
$d_{12}+d_{13}+d_{14}$.
An axis sum is the sum of labels of $2$ opposite edges --
e.g. $d_{12}+d_{34}$.

\begin{center}
\resizebox{!}{1.5in}{\includegraphics{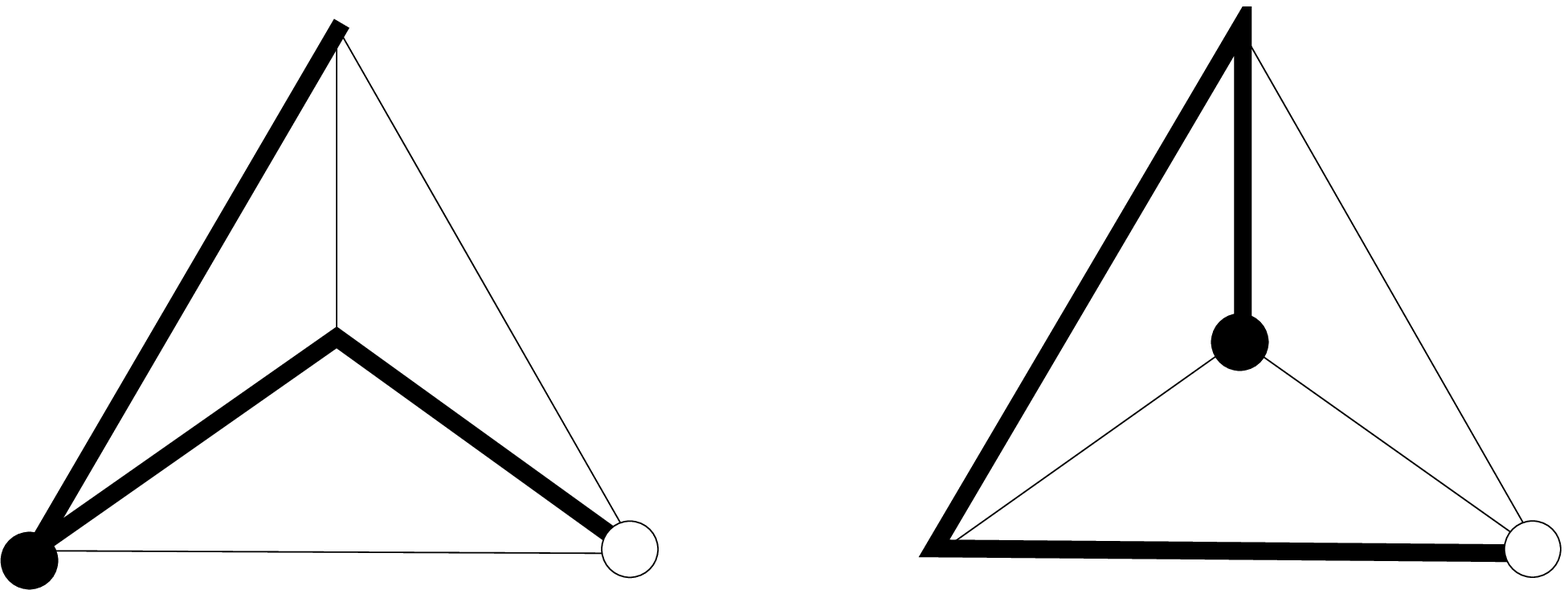}}
\newline
{\bf Figure 1.1:\/} $2$ of the $48$ decorations of $K_4$.
\end{center}

Figure 1.1 shows $2$ of the $48$ possible decorations
of $K_4$ in which we choose an embedded $3$-path, a
white endpoint of the path, and a black vertex of the
path which is not adjacent to the white endpoint.
For each such decoration $D$, we have a subset
$X_D \subset X$ consisting of those pseudo-tetrahedra
with the following properties:
\begin{itemize}
\item The axis sum of the opposite pair contained in $D$ is largest.
\item The axis sum of the opposite pair disjoint from $D$ is smallest.
\item The vertex sum at the black vertex is smallest.
\item The vertex sum at the white vertex is not greater than the
vertex sum at the vertex of $D$ incident to the white vertex.
\end{itemize}
It turns out that $X_D$ is linearly isomorphic to an orthant in $\R^6$. 
We call $X_D$ a
{\it chamber\/}.  Our construction partitions $X$ into
$48$ chambers.  We will explore this
partition more thoroughly in \S 2.

Let $\beta \subset K_4$ denote a subset of edges.
Call $\beta$ {\it friendly\/} if $K_4-\beta$
is not a union of edges all incident to the
same vertex. Otherwise, we call
$\beta$ {\it unfriendly\/}.  Up to isometry,
there are $7$ friendly subsets and $3$ unfriendly ones.

\begin{theorem}
\label{schema}
Let $\beta \subset K_4$ denote any friendly subset.
Let $g=D_{\beta}f$ denote the directional derivative of
$f$ along $\beta$. There is a nonempty union $X_{\beta}$ of
chambers of $X$, and constants $A_{\alpha}<B_{\beta}$,
with $B_{\beta}>0$, such that the function
$g \sum_{i<j} d_{ij} - C f$
is non-negative on $X_{\beta}$ if and only if
$C \in [A_{\alpha},B_{\beta}]$. Moreover,
every chamber of $X-X_{\beta}$ contains a point
where $f>0$ and $g<0$.
\end{theorem}

When $\beta=K_4$, Theorem \ref{ineq1} tells us that
$X_{\beta}=X$ and $(A_{\beta},B_{\beta})=(16,36)$.
Here is a summary of what we prove in the remaining cases.

\begin{itemize}

\item When $\beta$ is a single edge,
$X_{\beta}$ consists of the $12$
chambers $X_D$ such that $\beta \not \subset D$ and
the black vertex of $D$ is an endpoint of $\beta$.
(See Figure 5.1.) Here
$(A_{\beta},B_{\beta})=(0,2)$.

\item When $\beta$ is a pair of incident edges,
$X_{\beta}$ is the set of $4$ chambers $X_D$ such
that $\beta$ is disjoint from the outer two
edges of $D$ and the black dot is incident to
both edges of $\beta$. (See Figure 5.2) Here
$(A_{\beta},B_{\beta})=(0,12)$.

\item When $\beta$ is a pair of opposite edges,
$X_{\beta}$ is the set of $32$ chambers
$X_D$ such that $\beta \not \subset D$.
(See Figure 5.3 for some pictures.)
The constants satisfy $A_{\beta} \leq 0$ and $B_{\beta} \geq 4$.

\item When $\beta$ is $3$ edges incident to a
vertex $v$, the set 
$X_{\beta}$ consists of the $12$ chambers
$X_D$ such that the black vertex is $v$.
Here $(A_{\beta},B_{\beta})=(8,18)$.

\item When $\beta$ is a $3$-path, 
$X_{\beta}$ is the set of $8$ chambers
$X_D$ such that the black vertex is an
interior vertex of $\beta$, and
the outer edges of $D$ are disjoint from
$\beta$.  Here
$A_{\beta} \leq -6$ and $B_{\beta} \geq 16$.

\item When $\beta$ is a $4$-cycle, $X_{\beta}$ consists 
the $16$ chambers $X_D$ such that the outer two edges
of $D$ are disjoint from $\beta$. Here
$(A_{\beta},B_{\beta})=(0,24)$.

\end{itemize}

We have declared $3$-cycles unfriendly, but actually
we can say a lot about what happens for them.  When
$\beta$ is a $3$-cycle, let $X_{\beta}$ denote the $36$ chambers
$X_D$ so that the black vertex lies in $\beta$.

\begin{theorem}
\label{schema3}
Let $\beta$ be a $3$-cycle. 
Let $g=D_{\beta}f$ denote the directional derivative of
$f$ along $\beta$. Then the function
$g \sum_{i<j} d_{ij} - C f$
is non-negative on $X_{\beta}$ if and only if $C=8$.
Moreover, every chamber of $X-X_{\beta}$ contains a point
where $f>0$ and $g<0$.
\end{theorem}

\noindent
{\bf Remark:\/}
Roughly speaking, the decorations
defining $X_{\beta}$ try as hard as possible to
have their edges disjoint from $\beta$, and
their marked vertices contained in $\beta$.
\newline

Say that a lengthening
of a tetrahedron along a subset of edges
{\it locally increases\/} (respectively {\it decreases\/})
 the volume if the volume
goes up (respectively down) when we add the same sufficiently small
amount to each edge in the subset. 
Theorems \ref{schema} and \ref{schema3} have
the following immediate corollary.

\begin{corollary}
\label{schema2}
Suppose that $\beta$ is either a $3$-cycle or
a friendly subset of $K_4$.  
For any tetrahedron in $X_{\beta}$, the lengthening
along $\beta$ locally increases volume.  Moreover,
every chamber of $X-X_{\beta}$ contains 
a tetrahedron such that lengthening along
$\beta$ decreases volume.
\end{corollary}

One could say that Corollary \ref{schema2} gives coarsely
sharp conditions on when selective lengthening increases
volume.  Of course, if we used a finer triangulation,
we could make finer statements about this.

We can get some weak partial results about the remaining
two unfriendly configurations, those whose complement
is either a single edge or a pair of incident edges.
We will discuss this  briefly at
the end of \S 5.  Our methods really do fail for
these two configurations.
\newline
\newline
This paper has a companion computer program -- a
heavily documented and open-sourced
graphical user interface -- which the
interested readers can download from
\newline
\indent\indent
{\bf http://www.math.brown.edu/$\sim$res/Java/CM2.tar\/}
\newline
The program does all the integer polynomial
calculations, and also
shows plots of the Cayley-Menger determinant
and the various relevant directional derivatives, 
I discovered essentially everything 
in the paper using the program.
\newline
\newline
Here is some speculation on related questions.
Genevieve Walsh asked about results similar to
Theorems \ref{main}  and \ref{exist} 
for other combinatorial types - e.g. the octahedron.
One sensible constraint is that the maximum valence should be 
at most $5$, so that the combinatorial type
can be realized as convex polyhedra with equilateral
facets.   There are explicit analogues of the
Cayley-Menger determinant, which give volume
formulas for other combinatorial types. See [{\bf Sa\/}].

One might also ask about hyperbolic geometry versions 
of the results here. It seems
that Schlafli's formula -- see e.g. [{\bf L\/}] -- might be
useful.  It would be very nice to prove hyperbolic
or spherical versions of these results, and then deduce
the Euclidean results as limiting cases.  I have no idea
how to do this.

Just as Theorem \ref{exist} is a generalization of
Theorem \ref{main}, I wonder if Theorem \ref{schema}
and Corollary \ref{schema2} have
higher dimensional generalizations. It would be nice
to find a conceptual proof of Theorem \ref{schema}, because
my techniques are unlikely to be feasible in higher
dimensions.
\newline

I learned about Daryl Cooper's question during a lively semester
program in computational geometry, topology, and dynamics
at ICERM in Fall 2013.
I thank Bob Connelly, Peter Doyle, Ramin Naimi, Igor Rivin, Sinai Robins, and
Genevieve Walsh for interesting and helpful conversations
about this problem, some at ICERM and some elsewhere.  I would
especially like to acknowledge some conversations with 
Peter Doyle which helped guide me towards the special
$48$-chamber decomposition of the space $X$.  Peter made
the great guess that the vertex and axis sums should be
important in this edge-lengthening business.

\newpage

%% file: 2setup.tex
\section{Pseudo Tetrahedra}

\subsection{Normalized Pseudo Tetrahedra}

We say that a pseudo-tetrahedron is
{\it normalized\/} if
\begin{equation}
\sum_{i<j} d_{ij}=24.
\end{equation}
Let $X_{24}$ denote the space of normalized
pseudo-tetrahedra.
The $48$-partition of $X$ discussed in the 
introduction is the cone over a partition of
$X_{24}$ into $48$ $5$-simplices.  We choose
the normalization $24$ because it is the smallest
number we can choose which makes all these 
simplices integral.
Since all the inequalities we stated in
the introduction are homogeneous, it suffices to
prove them on $X_{24}$.

There are $7$ special points of $X_{24}$:
\begin{itemize}
\item $3$ of these points correspond to degenerate
tetrahedra in which the points have collapsed in pairs.
\item $4$ of these points correspond to degenerate
tetrahedra in which $3$ of the points have collapsed to one.
\end{itemize}
The $7$ vectors corresponding to these points are
\begin{itemize}
\item $A_1=(0,6,6,6,6,0)$.
\item $A_2=(6,0,6,6,0,6)$.
\item $A_3=(6,6,0,0,6,6)$.
\item $B_1=(8,8,8,0,0,0)$.
\item $B_2=(8,0,0,8,8,0)$.
\item $B_3=(0,8,0,8,0,8)$.
\item $B_4=(0,0,8,0,8,8)$.
\end{itemize}
We call these points {\it extrema\/} of
$X_{24}$, for reasons which will become clear
momentarily.

The following result is somewhat surprising, because the
points above are all (degenerate) tetrahedra whereas
$X_{24}$ certainly contains pseudo-tetrahedra which
are not tetrahedra in any sense -- e.g.
$(6,3,3,3,3,6)$.

\begin{lemma}
\label{hull}
$X_{24}$ is the convex hull of the $7$ extrema.
\end{lemma}

\startproof
Let $C$ denote the convex hull of the extrema.
Certainly $C \subset X$.  
For $k=1,2,3$, let
$\A_k$ denote the convex hull of the list of
$6$ extrema obtained by omitting $A_k$.
Each $\A_k$ is a $5$-simplex, and
$\A_i \cap \A_j$ is a $4$ simplex.
It is easy to check that $\A_i$ and $\A_j$
lie on opposite sides of the $4$-plane containing
their intersection.  

Call a face of $\A_j$ {\it free\/} if it is not also
a face of $\A_j$ for $j \not = i$. Otherwise, call
the face {\it bound\/}.  We have already exhibited
$2$ bound faces of each $\A_j$.  For the remaining
faces, we check that the barycenter of the face
lies in $\partial X_{24}$.  For instance, one of
the barycenters of a free face of $\A_1$ is
$$(A_2+A_3+B_1+B_2+B_3)/5=(28,22,14,22,14,20)/5.$$
In particular $d_{14}+d_{24}=d_{12}$.  The barycenter
condition implies that the entire free face lies in
$\partial X$.  This
$\bigcup \A_j$ is a union of three $5$-simplies,
with pairwise disjoint interiors, whose boundary
lies in $\partial X$.  This is only possible if
$X=\bigcup \A_j$.  But $\bigcup \A_j \subset C$.
Hence $C=X$.
\endproof

The proof in the Lemma \ref{hull} shows that
$X_{24}$ has a partition into $3$ simplices.
Let us consider the structure of $\A_1$.
If we compute the axis sums of the labelings
corresponding to the extrema, we find that
these sums are all equal for the $B$-extrema,
and $(12)(34)$ has largest axis sum for
$A_2$ and $A_3$.  Thus $\A_1$ consists
entirely of points whose largest axis sum
is $(12)(34)$.  Similarly,
$\A_2$ consists entirely of points whose
largest axis sum is $(13)(24)$, and
$\A_3$ consists entirely of points whose
largest axis sum is $(14)(23)$.

We can also define the simplex $\B_k$, which
is the convex hull of the list of $6$ extrema
obtained by omitting $B_k$.  The same proof
as above show that this gives a $4$-partition
of $X_{24}$.  An analysis similar to what
we did for the $3$-partition shows that
$\B_k$ consists of those points whose
corresponding labelings of $K_4$ have
smallest vertex sum at vertex $k$.
\newline
\newline
{\bf Remark:\/} There is a beautiful lower-dimensional
picture which gives a good feel for how the $3$-partition
and the $4$-partition are related.  One can think of
a triangular bi-pyramid $T$ as the join of a triangle and
a pair of points.  Correspondingly, $T$ has a partition
into $2$ tetrahedra, and also a partition into $3$
tetrahedra. This is the famous $2-3$ relation often
discussed in connection with $3$-dimensional triangulations.
The situation we have is a higher dimensional
analogue of this.

\subsection{The Common Refinement}

The space $X_{24}$ has a partition into $12$ simplices,
as follows:  We define $\C_{ij}$ to be the convex hull
of the point
\begin{equation}
C=(4,4,4,4,4,4)
\end{equation}
and the list of $5$ extrema obtained by omitting
$A_i$ and $B_j$.  

\begin{lemma}
\label{hull2}
$\C_{ij} = \A_i \cap \B_j$.
\end{lemma}

\startproof
Note that 
$$C=\frac{1}{3} \sum A_i = \frac{1}{4} \sum B_j,$$
So that $C \subset \A_i$ and $C \subset \B_j$
for all $j$.  Hence, all vertices of
$\C_{ij}$ are contained in $\A_i \cap \B_j$.
Hence $\C_{ij} \subset \A_i \cap \B_j$.

Next, we check that the barycenter of every face
of $\A_{ij}$ lies in $\partial(\A_i \cap \B_j)$.
By symmetry, it suffices to 
check this for $(i,j)=(1,4)$.  Again by symmetry,
it suffices to make the check for the face of
which does not involve $C$ and for one additional
face.  The barycenter of the face not involving
$C$ is the same as the one we computed in
the proof of Lemma \ref{hull}.  This point
must lie in both $\partial \A_i$ and
$\partial \B_4$ because it lies in 
$\partial X_{24}$.  One of the other barycenters is
$$(A_2+A_3+B_1+B_2+C)/5=(32,18,18,18,18,16)/5.$$
This point lies in $\partial \B_4$ because the
vertex sum at vertex $4$ is the same as the
vertex sum at vertex $3$.  This check establishes
what we want.

As in the proof of Lemma \ref{hull}, the barycenter
condition implies that every face of
$\C_{ij}$ lies in $\partial (\A_i \cap \B_j)$.
Since both sets are $5$-dimensional convex
polytopes, this situation is only possible if
$\C_{ij}=\A_i \cap \B_j$.
\endproof

Lemma \ref{hull2} implies that
$X_{24}$ has a partition into $12$ simplices,
namely $\C_{ij}$ for $i \in \{1,2,3\}$ and
$j \in \{1,2,3,4\}$.  The simplex
$\C_{ij}$ consists of those labelings where
the $i$th vertex sum is smallest and the
$j$th axis sum is largest.
\newline
\newline
{\bf Remark:\/} 
Once again, the picture for the bi-pyramid is
useful here.  The intersections of the $2$-partition
of the bi-pyramid with the $3$-partition gives a
$6$-partition into smaller tetrahedra.  Our situation
here is a higher dimensional analogue of this.

\subsection{The Final Partition}

The order $24$ symmetric group $S_4$ acts on the
space $X_{24}$ {\it via\/} label permutation.
The even subgroup $A_{4}$ acts freely and
transitively on our $12$-partition.   However,
the full group $S_5$ does not act freely on
the $12$-partition.  The stabilizer of each
simplex is an order $2$ subgroup.
We will use this symmetry to facilitate the understanding
of a refinement of the $12$-partition into a $48$-partition.
Basically, we cut each of the simplices into $4$ symmetric
pieces, again simplices.  We do this for the simplex
$\C_{11}$ and then use the $A_4$ symmetry to do it
for the remaining simplices.

We let $A_{ij}=(A_i+A_j)/2$ and likewise 
$B_{ij}=(B_i+B_j)/2$.  Also, we let $H(\cdot)$ stand
the for convex hull.
We introduce the $4$ simplices

\begin{equation}
\D_{1111}={\rm H\/}(C,B_2,B_{34},A_{23},B_3,A_2).
\end{equation}

\begin{equation}
\D_{1112}={\rm H\/}(C,B_2,B_{34},A_{23},B_3,A_3).
\end{equation}

\begin{equation}
\D_{1121}={\rm H\/}(C,B_2,B_{34},A_{23},B_4,A_2).
\end{equation}

\begin{equation}
\D_{1121}={\rm H\/}(C,B_2,B_{34},A_{23},B_4,A_3).
\end{equation}

Notice that only the last two vectors are changing.
One can see direcly that $\bigcup \D_{11ij}$ gives
a partition of $\D_{11}$.  What we are doing is
subdividing the $3$-simplex $H(B_3,B_4,A_2,A_3)$ into
$4$ symmetric pieces, and then taking the join with
the segment $H(C,B_2)$.  Figure 2.1 shows how to think
about the subdivision of the tetrahedron.

\begin{center}
\resizebox{!}{2in}{\includegraphics{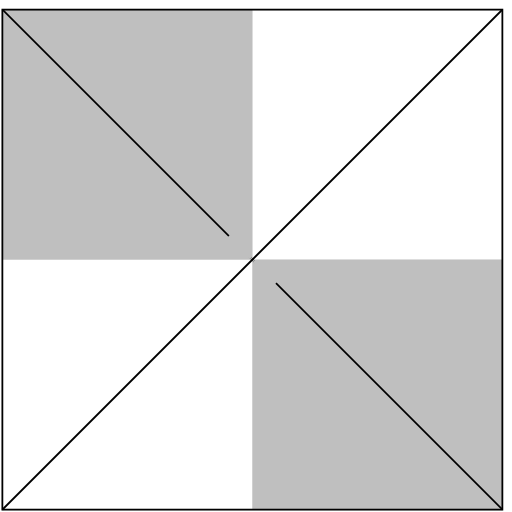}}
\newline
{\bf Figure 2.1:\/} Top view of the $4$-subdivision of a tetrahedron.
\end{center}

As we mentioned above, we use the $A_4$ symmetry to promote
our partition of $\D_{11}$ into a partition of all of
$X_{24}$ into $48$ simplices.

A direct calculation shows that the vertices of
$\D_{11ij}$ all satisfiy the inequalities associated
to the decorations discussed in the introduction.
Figure 2.2 shows the $4$ decorations.

\begin{center}
\resizebox{!}{4.5in}{\includegraphics{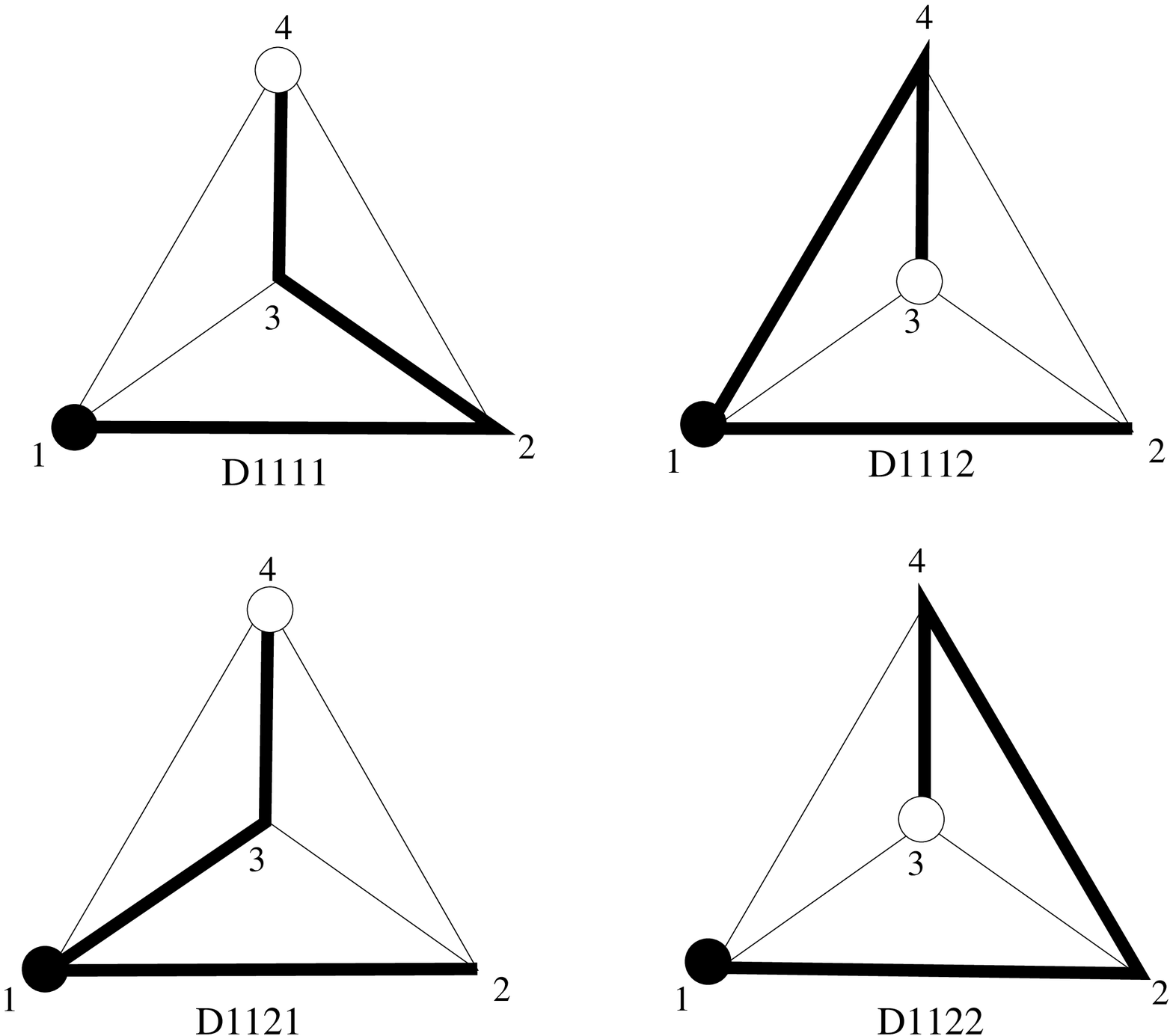}}
\newline
{\bf Figure 2.2:\/} The decorations associated to
$\D_{11ij}$.
\end{center}

The partition discussed in the introduction restricts
to a partition of $X_{24}$ into $48$ convex polytopes.
Each of these polytopes contains one of our sets
$\D_{ijk\ell}$.  But then the partition in the
introduction must intersect $X_{48}$ precisely in
the partition $\{\D_{ijk\ell}\}$ constructed here.
By honogentity, each chamber in the partition
from the introduction is the cone over some
simplex in our partition here.  This establishes
the claim in the introduction that the chambers
are linearly isomorphic to orthants.

\subsection{Theorem \ref{ineq1} Modulo a Detail}

Let $f$ denote the Cayley-Menger determinant, as in
Equation \ref{cm}, and let $g$ be as in
Equation \ref{df}.  In \S 3 we will explain how
we establish the following theorem:

\begin{lemma}
\label{norm}
$2g-3f$ and $3g-2f$ are non-negative on $\D_{11}$.
\end{lemma}

\begin{corollary}
$2g-3f$ and $3g-2f$ are non-negative on $X_{24}$.
\end{corollary}

\startproof
The two functions $f$ and $g$ are invariant under
the action of the alternating group $A_4$, and
this action freely permutes the simplices in
the $12$-partition of $X_{24}$.  Since we have
non-negativity on one of these simplices, we get
the non-negativity on all of them, by symmetry.
\endproof

\begin{lemma}
\label{confiner}
Let $a$ and $b$ be constants.  Then
$bg+af \geq 0$ on $X_{24}$ if and only if
$bg+af$ is a non-negative combination of the
two functions in Lemma \ref{norm}.
\end{lemma}

\startproof
Recall that $C=(4,4,4,4,4,4)$.
We compute that
\begin{equation}
\label{concrete1}
f(C)=3^0 \times 2^{14}, \hskip 30 pt
g(C)=3^1 \times 2^{13}.
\end{equation}
This equation shows that $2b+3a \geq 0$.

We compute that
\begin{equation}
\label{concrete2}
f(6,3,3,3,3,6)=-3^6 \times 2^{7}, \hskip 30 pt
g(6,3,3,3,3,6)=-3^5 \times 2^{8}.
\end{equation}
This equation shows that $-3b-2a \geq 0$.

Finally, by considering a degenerate tetrahedron
with $3$-fold symmetry and an equilateral triangle
base, we see that it can happen that $f=0$ and $g>0$.
This forces $b \geq 0$. Equation \ref{concrete2}
then force $a \leq 0$.  Our two inequalities above
now confine $(a,b)$ to one half of a cone in the
plane, and this cone is precisely the set 
of coeffients one can obtain by taking non-linear
combinations of the ones in Lemma \ref{norm}.
\endproof

Combining the results in this section, with some
basic arithmetic, and the fact that
$\sum_{i<j}d_{ij}=24$, we see that
Theorem \ref{ineq1} is true on $X_{24}$.
But then, by homogeneity, Theorem \ref{ineq1} is
true on all of $X$.

\subsection{Proof of Theorem \ref{main}}

To prove Theorem \ref{main}, we find it
convenient to work on the subset $X_6 \subset X$
consisting of pseudo-tetrahedra $\{d_{ij}\}$
where $\sum_{i<j} d_{ij}=6$.

\begin{lemma}
\label{scale1}
On $X_6$ we have
$g \geq 6f$.
\end{lemma}

\startproof
This is a special case of Theorem \ref{ineq1}.
\endproof

Rather than prove Theorem \ref{main} for the
unit lengthening of a tetrahedron, we will
prove that the $t$-lenghening $\Delta_t$ of
a tetrahedron $\Delta_0 \in X_6$ satisfies the
volume bound
\begin{equation}
\label{elaborate}
\frac{{\rm volume\/}(\Delta_t)}{{\rm volume\/}(\Delta_0)} 
\geq (1+t)^3.
\end{equation}
Theorem \ref{main} follows from this result, and scaling.

Let $\phi_t:  X \to  X$ denote
the flow defined by 
\begin{equation}
\phi_t(d_{ij})=\{d_{ij}+t\}.
\end{equation}
Let $D \in X_6$ represent $\Delta_0$.   Let
$D_t=\phi_t(D)$ and $F(t)=f(D_t)$.

Since $f$ is homogeneous of degree $6$ and $g$ is homogeneous
of degree $5$, Lemma \ref{scale1} implies that
\begin{equation}
\label{diffeq}
\frac{dF}{dt} \geq \frac{6}{1+t}F(t).
\end{equation}
This shows immediately that $F$ is increasing.
Suppose that $D_t$ fails to be tetrahedral for some $t>0$.
This would force $F(t)=0$, contradicting the increase of $F$.
This shows that $D_t$ is tetrahedral for all $t>0$.

Equation \ref{diffeq} can be rearranged as
\begin{equation}
\label{compare}
\frac{d}{dt} \log F \geq \frac{6}{1+t}.
\end{equation}
Integrating both sides and then exponentiating, we get
\begin{equation}
\frac{F(t)}{F(0)} \geq (1+t)^6.
\end{equation}
Taking square roots and using Equation \ref{classic},
we get exactly Equation \ref{elaborate}.

\newpage

%% file: 3posdom.tex
\section{The Method of Positive Dominance}

In \S 2 we reduced Theorem \ref{main} to
Theorem \ref{ineq1}.  In this chapter we will
explain our computational method for proving
Theorem \ref{ineq1}.  The material in this
chapter is taken mostly from my
recent monograph [{\bf S\/}], though
it has been adapted to the present situation.

\subsection{Single Variable Case}

As a warmup, we consider the situation for polynomials
in a single variable.
Let 
\begin{equation}
\label{poly}
P(x)=a_0+a_1x+...+a_nx^n
\end{equation}
be a polynomial with real coefficients.  Here we describe
a method for showing that $P \geq 0$ on $[0,1]$,

Define
\begin{equation}
\label{sum1}
A_k=a_0+ \cdots + a_k.
\end{equation}
We call $P$ {\it weak positive dominant\/} (or
{\it WPD\/} for short) if
$A_k \geq 0$ for all $k$.
\newline
\newline
{\bf Remark:\/}
To keep consistent with [{\bf S\/}], we 
reserve the terminology {\it positive dominant\/}
for the case $A_k >0$ for all $k$.  However,
in this paper we only care about weak positive dominance.

\begin{lemma}
\label{PD1}
If $P$ is weak positive dominant, then
$P \geq 0$ on $[0,1]$.
\end{lemma}

\startproof
The proof goes by induction on the degree of $P$.
The case $\deg(P)=0$ follows from the fact that $a_0=A_0 \geq 0$.
Let $x \in [0,1]$.
We have
$$P(x)=a_0+a_1x+x_2x^2+ \cdots + a_nx^n \geq $$
$$a_0x+a_1x+a_2x^2+ \cdots + a_nx^n=$$
$$x(A_1+a_2x+a_3x^2+ \cdots a_nx^{n-1})=xQ(x) \geq 0$$
Here $Q(x)$ is weak positive dominant and has degree $n-1$.
\endproof

\noindent
{\bf Remark:\/}
The converse of Lemma \ref{PD1} is generally false.
We will give an example below.
\newline

Given an interval $I \subset \R$, let $A_I$ be the affine and
orientation preserving map which carries $[0,1]$ to $I$. We
call the pair $(P,I)$ {\it weak positive dominant\/} if
$P \circ A_I$ is WPD. If $(P,I)$ is WPD
then $P \geq 0$ on $I$, by Lemma \ref{PD1}.  

We say that a partition 
$[0,1]=I_1 \cup ... \cup I_n$
is {\it weak positive dominant with respect to \/} $P$ 
if $(P,I_k)$ is WPD for each $k=1,...,n$.
For short, we will just say that $P$ 
{\it has a weak positive dominant partition\/}.
If $P$ has a WPD partition, then
$P \geq 0$ on $[0,1]$.
\newline
\newline
{\bf Example:\/}
The polynomial $$P(x)=3-4x+2x^2$$ is not WPD
but satisfies $P(x) \geq 1$ for all $x \in \R$.
Consider the partition $[0,1]=I_1 \cup I_2$, where $I_1=[0,1/2]$
and $I_2=[1/2,1]$.  The corresponding affine maps are
$$A_1(x)=x/2; \hskip 30 pt
A_2(x)=x/2+1/2.$$
We compute
$$P \circ A_1(x)=3-2x+x^2/2, \hskip 30 pt
P \circ A_2(x)=3/2-x+x^2/2.$$
Both of these polynomials are WPD.
Hence $P$ has a WPD partition.
\newline
\newline
\noindent
{\bf Divide-and-Conquer Algorithm:\/}
If $P\geq 0$ on $[0,1]$, we can try to find a WPD
partition using a divide-and-conquer
algorithm.  The algorithm works like this.
\begin{enumerate}
\item Start with a list LIST of intervals. Initially
LIST consists only of $[0,1]$. 
\item Let $I$ be the last
interval on LIST. We delete $I$ from LIST and
then test
whether $(P,I)$ is weak positive dominant.
\item Suppose $(P,I)$ is weak positive dominant. 
We go back to Step 2 if LIST is nonempty
and otherwise halt.
\item Suppose $(P,I)$ is not weak positive dominant.
We append to LIST the two intervals
obtained from cutting $I$ in half, then
go back to Step 2.
\end{enumerate}
If the algorithm halts, then (assuming
that the calculations are done exactly)
we have a proof that $P \geq 0$ on $[0,1]$.

\subsection{The General Case}

Now we go to the higher dimensional case.
We consider real polynomials in the variables
$x_1,...,x_k$.
Given a multi-index 
$I=(i_1,...,i_k) \in (\N \cup \{0\})^k$ we
let 
\begin{equation}
x^I=x_1^{i_1}...x_k^{i_k}.
\end{equation}
Any polynomial $F \in \R[x_1,...,x_k]$
can be written succinctly as
\begin{equation}
F=\sum A_I X^I, \hskip 30 pt A_I \in \R.
\end{equation}
If $I'=(i_1',...,i_k')$ we write
$I' \leq I$ if $i'_j \leq i_j$ for
all $j=1,...,k$.
We call $F$ {\it weak positive dominant\/}
if
\begin{equation}
\label{summa}
\sum_{I' \leq I} A_{I'} \geq 0
\hskip 30 pt \forall I,
\end{equation}

\begin{lemma}
\label{PD2}
If $P$ is weak positive dominant then $P \geq 0$ on $[0,1]^k$.
\end{lemma}

\startproof
The $1$ variable case is
Lemma \ref{PD1}.
In general, we write
\begin{equation}
P=f_0+f_1x_k+...+f_mx_k^m,
\hskip 20 pt f_j \in \R[x_1,...,x_{k-1}].
\end{equation} 
Let $P_j=f_0+...+f_j$.
Since $P$ is weak positive dominant, we get that
$P_j$ is weak positive dominant for all $j$.
By induction on $k$, we get $P_j \geq 0$ on
$(0,1)^{k-1}$. 
But now, if we hold
$x_1,...,x_{k-1}$ fixed and let
$t=x_k$ vary, the polynomial
$g(t)=P(x_1,...,x_{k-1},t)$
is weak positive dominant..  Hence, by Lemma \ref{PD1},
we get $g \geq 0$ on $[0,1]$.
Hence $P \geq 0$
on $[0,1]^k$.
\endproof

We can perform the same kind of divide-and-conquer
algorithm as in the $1$-dimensional case.
We always take our domain to be $[0,1]^k$. 
Let $P$ be a polynomial.  We are going to describe
our subdivision in terms of what it does to the
polynomials rather than what it does to the domain.
 
We first define the maps
\begin{equation}
A_j(x_1,...,x_k)=\bigg(x_1,..x_{j-1}.,\frac{x_j}{2},x_{j+1}...,x_k\bigg).
\end{equation}

\begin{equation}
B_j(x_1,...,x_j,...,x_k)=(x_1,...,x_{j-1},1-x_j,x_{j+1}...,x_k).
\end{equation}

We define the $j$th {\it subdivision\/} of $P$ to be the
set
\begin{equation}
\label{SUB}
\{P_{j1},P_{j2}\}=
\{P \circ A_j,\ P \circ B_j \circ A_j\}.
\end{equation}

\begin{lemma}
\label{half}
$P \leq 0$ on $[0,1]^k$ if and only if
$P_{j1} \geq 0$ and $P_{j2} \geq 0$ on
$[0,1]^k$.
\end{lemma}

\startproof
By symmetry, it suffices to take $j=1$.
Define
\begin{equation}
[0,1]^k_1=[0,1/2] \times [0,1]^{k-1}, \hskip 30 pt
[0,1]^k_2=[1/2,1] \times [0,1]^{k-1}.
\end{equation}
Note that
\begin{equation}
A_1([0,1]^k)=[0,1]^k_1,\hskip 30 pt
B_1 \circ A_1([0,1]^k)=[0,1]^k_2.
\end{equation}
Therefore,
$P \geq 0$ on $[0,1]^k_1$ 
if and only if $P_{j1} \geq 0$ on $[0,1]^k$. Likewise
$P \geq 0$ on $[0,1]^k_2$ if and only if
if $P_{j2} \geq 0$ on $[0,1]^k$.
\endproof

Say that a {\it marker\/} is a non-negative
integer vector in $\R^k$. 
Say that the {\it youngest entry\/} in the the marker
is the first minimum entry going from left to right. 
The {\it successor\/} of a marker is the marker obtained
by adding one to the youngest entry. For instance,
the successor of $(2,2,1,1,1)$ is $(2,2,2,1,1)$.
Let $\mu_+$ denote the successor of $\mu$.

We say that a {\it marked polynomial\/} is a pair
$(P,\mu)$, where $P$ is a polynomial and
$\mu$ is a marker.  Let $j$ be the position of the
youngest entry of $\mu$.  We define the
{\it subdivision\/} of $(P,\mu)$ to be the
pair
\begin{equation}
\{(P_{j1},\mu_+,(P_{j2},\mu_-)\}.
\end{equation}
Geometrically, we are cutting the domain in half
along the longest side, and using a particular rule
to break ties when they occur.
\newline
\newline
{\bf Divide-and-Conquer Algorithm:\/}
\begin{enumerate}
\item Start with a list LIST of marked polynomials.
Initially, LIST consists only of the marked polynomial
$(P,(0,...,0))$.
\item Let $(Q,\mu)$ be the last element of LIST.
We delete $(Q,\mu)$ from LIST and test whether
$Q$ is weak positive dominant.
\item Suppose $Q$ is weak positive dominant.
we go back to Step 2 if LIST is not empty.
Otherwise, we halt.
\item Suppose $Q$ is not weak positive dominant.
we append to LIST the two marked polynomials
in the subdivision of $(Q,\mu)$ and then go
to Step 2.
\end{enumerate}

We call $P$ {\it Recursively Weak Positive Dominant\/}
or (RWPD) if the divide and conquer algorithm
halts for $P$.  If $P$ is RWPD then $P \geq 0$
on $[0,1]^k$.  This is a consequence of
Lemma \ref{half} and induction on the number of
steps taken in the algorithm.

\subsection{From Cubes to Simplices}
\label{cubetransform}

So far we have been talking about showing that
polynomials are non-negative on the unit cube
$[0,1]^k$.  But, we really want to show that
the polynomials of interest to us, namely
those from Lemma \ref{norm}, are non-negative on
the simplex $\D_{11}$.
In this section,
we explain how this is done.  We set $k=5$ and
use coordinates $(a,b,c,d,e)$ on $\R^5$.
\newline
\newline
{\bf Cube To Standard Simplex:\/}
Let $S_5 \subset \R^5$ denote the simplex
\begin{equation}
\{(a,b,c,d,e)|\ 1 \geq a \geq b \geq c \geq d \geq e \geq 0\}.
\end{equation}
We call $S_5$ the {\it standard simplex\/}, though
actually we won't use this terminology after this section.
There is a polynomial surjective map from
$[0,1]^5$ to $S_5$:
\begin{equation}
U(a,b,c,d,e)=(a,ab,abc,abcd,abcde).
\end{equation}

\noindent
{\bf Standard Simplex to Regular Simplex:\/}
Let $\Delta_6 \subset \R^6$ denote the regular $5$-simplex
in $\R^6$ consisting of the convex hull of the
standard basis vectors. That is,
$\Delta_6$ consists of points
$(x_1,...,x_6)$ such that $x_j \geq 0$ for all $j$
and $\sum x_j=1$.
There is an affine isomorphism from
$S_5$ to $\Delta_6$:
\begin{equation}
\label{trans1}
V(a,b,c,d,e)=(1-a,a-b,b-c,c-d,d-e,e).
\end{equation}
The easiest way to see that this works is to check
it on the vertices of $S_5$.
\newline
\newline
\noindent
{\bf Regular Simplex to General Simplex:\/}
Let $\Sigma$ denote a $5$-simplex in
$\R^6$.  We can think of
$\Sigma$ as a $6 \times 6$ matrix whose
$6$ columns are the vertices of $\Sigma$.
Call this matrix $W_{\Sigma}$.  The map
$W_{\Sigma}$ gives an affine isomorphism
from $\Delta_6$ to $\Sigma$.
\newline

Note that the composition
\begin{equation}
Z_{\Sigma}=W_{\Sigma} \circ V \circ U: [0,1]^5 \to \Sigma
\end{equation}
is a surjective rational map.
Given a polynomial $P$ and a simplex $\Sigma \subset \R^5$,
we define the new polynomial
\begin{equation}
\label{pullback}
P_{\Sigma}=P \circ Z_{\Sigma}.
\end{equation}

By construction, $P \geq 0$ on $\Sigma$ provided
that $P_{\Sigma}$ is RWPD.

\subsection{Proof of Lemma \ref{norm}}

Let $P$ and $Q$ be the two polynomials from
Lemma \ref{norm} and let $\Sigma=\D_{11}$,
the simplex from Lemma \ref{norm}. 

\begin{lemma}
\label{aux2}
$P_{\Sigma}$ and $Q_{\Sigma}$ are both
RWPD.
\end{lemma}

\startproof
We prove Lemma \ref{aux2} simply by coding all
the algebra in sight into a Java program and running it.
In the next chapter we discuss the 
implementations of our calculations.
For the function $P_{\Sigma}$ the algorithm
takes $7455$ steps and runs in about $62$ minutes
on my $2012$ Macbook pro.  For the function
$Q_{\Sigma}$, the algorithm takes $1173$ steps
and runs in about $6$ minutes.
\endproof

This proves Lemma \ref{norm}.  
\newline
\newline
{\bf Remarks:\/} For a given function,
the number of steps in the algorithm would be
the same on any perfectly running computer,
but of course the time would vary.  The number
of steps looks large, but one has to remember
that we are running the algorithm on a $5$ dimensional
cube.  Note that $(2^3)^5=32768$, so we are making
an average of less than $3$ subdivisions in each
coordinate direction.

\newpage

%% file: 4code.tex
\section{Implementation}

We implement our calculations in Java.  Here
we describe the salient features of the
code.

\subsection{Formulas for the Main Functions}
\label{form}

The reader can find explicit formulas for
the Cayley-Menger determinant $f$ and
the partial derivatives $\partial_j f$ for
$j=1,2,3,4,5,6$ in the file {\bf DataCM.java\/}.
This is one of the files in the directory you get
when you download my program from
\newline
\newline
\indent \indent \indent
{\bf http://www.math.brown.edu/$\sim$res/Java/CM2.tar\/}
\newline
\newline
Aside from the formulas for $f$ and its partial derivatives,
which we derived using Mathematica [{\bf W\/}], the rest of
the program is self-contained, in that all the
calculations for the paper are done there.  In an earlier
version of the paper and program, I implemented some of the
calculations in Mathematica, but that is no longer the case.

\subsection{General Features}

The program is written entirely in Java.  It has several
useful features.
\begin{itemize}
\item The user can see all the code.
\item The user can see all the decorations corresponding
to the $48$-simplex partition.
\item The user can run his/her own experiments, testing
various linear combinations of $f$ and its partial derivatives
for positivity.
\item While running, the program
has a documentation feature, so that the
user can learn about practically every facet of the program by
reading the text.
\item The program has a debugging tool which allows one to see that
the critical formulas really are correct.  
\end{itemize}
Aside from the debugging and the experiment modes, all the calculations having
to do with the proofs are done with exact integer arithmetic.

\subsection{Data Structures}

Here we describe the special data structures used
by our program.
\newline
\newline
\noindent
{\bf BigIntegers:\/} The BigInteger class in
Java is designed to do arbitrary digit
arithmetic.  In practice this means that
one can do arithmetic with integers which
have thousands of digits.  In our case,
we never get integers with more than, say,
$20$ digits.  So, we are well inside the
working power of the language.
\newline
\newline
{\bf Monomials:\/}
For us, a {\it MonoN\/} is a tuple $[C,(e_1,...,e_n)]$.
Here $C$ is a BigInteger and $e_1,...,e_n$ are non-negative
integers.   
This expression has the following meaning.
\begin{equation}
\label{monomial}
[C,(e_1,...,e_n)]=
C \prod_{i=1}^5 x_n^{e_n}.
\end{equation}
We only implement this class for $N=5$ and $N=6$,
and we found it convenient to define separate classes for
each case.
\newline
\newline
{\bf Polynomials:\/}
A {\it polyN\/} is a finite list of
monomials.  In practice we allow for
$50000$ (Mono6)s and $20000$ (Mono5)s.
For a Poly5, we have a universal
degree bound $e_j \leq 6$ for all $j$.
We don't keep track of the degree bound
for a Poly6, but the universal bound is
around $6$ as well.

\subsection{Poly6 Arithmetic}

We use the Poly6 class to implement the maps
discussed in \S \ref{cubetransform}.
The functions $f$ and $\partial_i f$ are stored as lists of
integers which the program readily converts into (Poly6)s.
Starting with a Poly6 P, which is some integer combination
of $f$ and its partial derivatives, and an integer
simplex $\Sigma\subset X_{24}$, we compute the
polynomial $P_{\Sigma}$.
The implementation is straightforward, and basically
is built out of polynomial addition and multiplication.
The routines are contained in the 
fairly well documented file {\bf Poly6.java\/}.

In principle, the function $P_{\Sigma}$ could
be precomputed using Mathematica, since we just
need to load in the formula once, and then process
it. Indeed, an earlier version of the program
did this.  However, since we wanted to test many
functions, we didn't want to repeatedly go into
Mathematica and save the output to our Java files.
My point is that the critical polynmomial
algebra we do is the implementation of the
Positive Dominance Algorithm.  

\subsection{Poly5 Arithmetic}

Here we describe the basic operations
we perform on polynomials.  The main point
is to implement the Positive Dominance Algorithm.
\newline
\newline
\noindent
{\bf Rotation:\/}
The $1$-{\it rotation\/} of
$[C,(e_1,...,e_5)]$
is $[C,(e_5,e_1,e_2,e_3,e_4)].$
The $k$-{\it rotation\/} is
obtained by applying the $1$-rotation
$k$ times.  The $k$-rotation
of a polyomial is simply the
list of $k$-rotations of its
monomials.  The $k$-rotation $R_kP$
of a polynomial $P$ is simply the
composition of $P$ with some
cyclic permutation of the coordinates.

If we have some operation $Z$ which
does something to the first coordinates
of the monomials of $P$, the
operation $R_{k-1}ZR_{-k+1}$ does the
same thing to the $k$th coordinate.
We use this trick so that we just
have to implement our main routines
for the first coordinate.  
\newline
\newline
{\bf Dilation:\/}
We call the operation in Equation
\ref{SUB} {\it dilation\/}.
We only implement dilation for the first variable.
Given a polynomial $P$, the two
polynomials from Equation \ref{SUB} are not
necessarily integer polynomials.  Their
coefficients are dyadic rationals: Expressions
of the form $p/2^{e_1}$. We let $E=\max e_j \leq 6$.
Where the expression is taken over all monomials.
To get an integer polynomial, we
use $2^E P_{11}$ and $2^E P_{12}$ instead
of $P_{11}$ and $P_{12}$.  
\newline
\newline
{\bf Reflection:\/}
Here {\it reflection\/} is the operation
of replacing the polynomial $P$ with the
polynomial $P \circ B_j$.  When $j$,
the new polynomial is
$P(1-a,b,c,d,e)$.  We only implement
the reflection operation for the first variable.
The reflection operation is  the
rate limiting step in our program, so we
explain the implementation carefully.
We first create a
$7 \times 7 \times 7 \times 7 \times 7$
array $\beta$ of BigIntegers. For each
term $[C,(e_1,...,e_5)]$ of the polynomial
we perform the following:  We let
$\tau$ be the $(e_1)$st row of
Pascal's triangle, with the sign
switched on the even terms.
For instance, rows $0$ and $1$ and $2$ are
$(0)$ and $(-1,1)$ and $(1,-2,-1)$ respectively.
We then make the substitution
\begin{equation}
\beta[j,e_2,e_3,e_4,e_5]=
\beta[j,e_2,e_3,e_4,e_5]+C \tau(j), \hskip 30 pt
j=0,...,e_1.
\end{equation}
Again, we always have $e_1 \in \{0,...,6\}$.

When we are done, we convert our array $\beta$ back
into a polynomial by including the monomial
$[C,(e_1,...,e_5)]$ iff the
final value of $\beta(e_1,...,e_5)$ is $C$.
The resulting polynomial is the reflection of $P$.
\newline
\newline
{\bf Subdivision:\/}
The polynomial $P_{k1}$ is obtained by the
following operations:
\begin{enumerate}
\item Let $Q$ be the $(-k)$th rotation of $P$.
\item Let $R$ be the dilation of $Q$.
\item Let $S$ be the $(k)$th rotation of $R$.
\item Return $S$.
\end{enumerate}
The polynomial $P_{k2}$ is obtained by the
following operations:
\begin{enumerate}
\item Let $Q$ be the $(-k)$th rotation of $P$.
\item Let $Q^*$ be the reflection of $Q$.
\item Let $R$ be the dilation of $Q^*$.
\item Let $S$ be the $(k)$th rotation of $R$.
\item Return $S$.
\end{enumerate}
In the case of $P_{12}$ it is very important
that steps 2 and 3 are not interchanged.

\subsection{Test for Weak Positive Dominance}

Building on the notation of Equation \ref{monomial}, we write
$$[C,(e_1,...,e_5)] \preceq (i_1,...,i_5)$$
if and only if $e_j \leq i_j$ for all $j$.
Let $E$ be the smallest multi-index so that
$[C,(e_1,...,e_5)] \preceq E$ for all terms.
We do the following loop.
\newline
\newline
for$(i_1=0;i_1\leq E_1;++i_1)\ \{$ \newline
for$(i_2=0;i_2\leq E_2;++i_2)\ \{$ \newline
for$(i_3=0;i_3\leq E_3;++i_3)\ \{$ \newline
for$(i_4=0;i_4\leq E_4;++i_4)\ \{$\newline
for$(i_5=0;i_5\leq E_5;++i_5)\ \{$ \newline
\indent Sum the coefficients
of all terms $[(c,(e_1,...,e_5)] \preceq (i_1,...,i_5)$. \newline
$\}\}\}\}\}$. \newline \newline
We return ``false'' if we ever get a negative total sum.
Otherwise we return ``true'', indicating that
the polynomial is weak positive dominant.
This procedure is probably rather far from being optimal,
but it is quite simple.

\subsection{Test for Negativity}

At each step of the WPDA, we perform an additional test.
We check whether or not the current polynomial is
negative at the origin.  If the current polynomial
is negative at the origin, we terminate the algorithm
because we have a proof that the original polynomial
is {\bf not\/} non-negative on the given simplex. 
In this way, our algorithm typically halts either with
a proof of non-negativity or a proof that some negative
values exist.

There is the theoretical possibility that we could
encounter a non-negative function that is not
recursively weak positive dominant.  For instance,
the polynomial $P(x,y)=(x-y)^2$ is non-negative on
the unit cube in any dimenson greater than $2$, but
not RWPD.  Fortunately, we do not encounter polynomials
like this in practice.  So, in all cases, our
algorithm terminates with a definite conclusion.

\subsection{Anti-Certification}
\label{anti}

Supposing that $\beta$ is some subset of edges, our
results also make statements about the chambers
of the space $X$ which do not belong to the union
$X_{\beta}$.   Let $Y$ be such a chamber.  On $Y$ we
want to show that it can happen that $f>0$ and $df<0$.
here $df$ is the directional derivative $D_{\beta}f$.
In this section, we explain how we do this rigorously.

We sample random points in $Y_{24}$ until we find a
candidate point $p \in Y_{24}$ such that (according
to floating point calculations, it appears that)
$f(p)>0$ and $df(p)<0$. 
We replace $p$ by a point $p^* \in Y \cap \Z^6$ in such a
way that $p$ and $p^*$ nearly lie on the same line
through the origin.  We then show that
$f(p*)>0$ and $df(p^*)<0$ using exact integer calculations.

It doesn't matter how we produce $p^*$,and it also doesn't
matter that $p$ and $p^*$ nearly lie on the same line
through the origin, but this property makes it likely
that $f(p^*)>0$ and $df(p^*)<0$.  Even though these
details don't matter from a logical standpoint, it
seems worth explaining how we get $p^*$. The point $p$ has
the form $L(q)$, where $L$ is an integer linear map
taking the standard simplex $\sum x_i=1$ to $Y_{24}$.
We then replace $q$ by the point
\begin{equation}
q^*={\rm floor\/}(10^{10}q)
\end{equation}
and set $p^*=L(q^*)$.
This does the job for us.

\newpage

%% file: 5select.tex
\section{Selective Lengthening}

\subsection{A Single Edge}
\label{singleedge}

The goal of this chapter is to prove
Theorem \ref{schema}, and to establish all
the supplementary facts mentioned after
we stated Theorem \ref{schema}.
We just go through the cases one at a time.

Let $\beta$ denote a single edge.  
Let $X_{\beta}$ be the union of $12$
chambers $X_D$ such that $\beta \not \subset D$ and
the black vertex of $D$ is an endpoint of $\beta$.
Figure 5.1 shows $3$ representative examples,
corresponding to the simplices listed below.
dThe edge $\beta$ is drawn in grey.  

\begin{center}
\resizebox{!}{1.4in}{\includegraphics{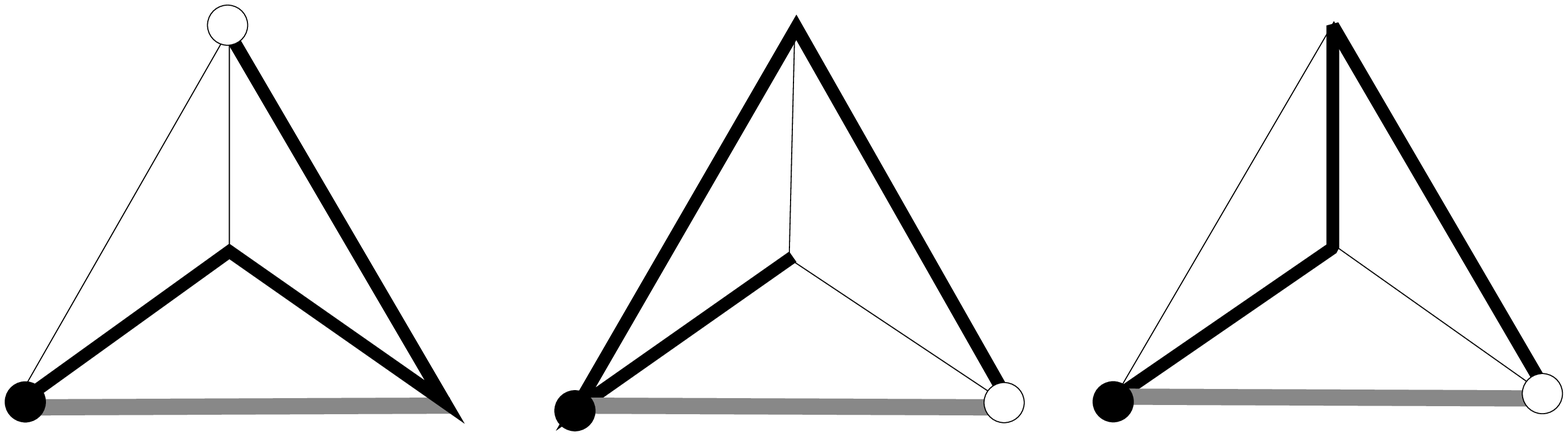}}
\newline
{\bf Figure 5.1:\/} Three of the $12$ decorations for $X_{\beta}$.
\end{center}

First of all, we use the anti-certification algorithm
discussed in \S \ref{anti} to show that any
chamber of $X-X_{\beta}$ has a point where
$f>0$ and $g<0$.  To save words, we will say below
that we {\it anti-certify the chambers not in $X_{\beta}$\/}.
 
Now we turn to $X_{\beta}$.
The $12$ chambers intersect the normalized space $X_{24}$ in
$12$ simplices.
If we pick $\beta=\{e_{12}\}$, then $3$ of these simplices
are given by
$$(4,4,4,4,4,4),\ (0,8,0,8,0,8),\ (4,0,4,4,8,4)$$
$$(3,6,3,3,6,3),\ (8,0,0,8,8,0),\ (0,6,6,6,6,0)$$
and
$$(4,4,4,4,4,4),\ (0,8,0,8,0,8),\ (4,0,4,4,8,4)$$
$$(3,6,3,3,6,3),\ (0,0,8,0,8,8),\ (0,6,6,6,6,0)$$
and
$$(4,4,4,4,4,4),\ (0,8,0,8,0,8),\ (4,0,4,4,8,4)$$
$$(3,6,3,3,6,3),\ (0,0,8,0,8,8),\ (6,6,0,0,6,6)$$
The remaining simplices are images of these under
the action of the subgroup of $S_4$ which
stabilizes $e_{12}$.

We work with the functions
\begin{equation}
P=g, \hskip 20 pt
Q=12 g - f.
\end{equation}
Again $g=D_{\beta}f$.  This function
depends on $\beta$, of course.
These are multiples of the ones mentioned above.
We use the Method of Positive Dominance to check
that both these functions are non-negative on
the $3$ simplices above.  By symmetry, $P$ and $Q$
are non-negative on the intersection of $X_{24}$
with $X_{\beta}$.  But then, by homogeneity, they
are non-negative on $X_{\beta}$.
For $P$, the number of steps taken for the
three simplices is $421,421,427$ respectively.
For $Q$, the number of steps taken for the
three simplices is $457,469,617$ respectively.

Our calculations show that the cone
$y=0$ and $y>x/12$ contains all points of
the form $(f,g)$ when this pair of functions
is evaluated on $X_{\beta}$. To show that
$(A_{\beta},B_{\beta})=(0,2)$, we just need to see that
there is no smaller cone which has this property.
In other words, if we tilt the two lines bounding
the cone inward, so to speak, the lines will cross
points of the image.  We deal with the two lines
in turn.

One of the edges of $X_{\beta}$ has endpoints
$(0,6,6,6,6,0)$ and $(3,6,3,3,6,3)$.  Consider
the point
\begin{equation}
\Omega_t=(1-t)(0,6,6,6,6,0)+t(3,6,3,3,6,3).
\end{equation}
We compute that
\begin{equation}
g(\Omega_t)+t f(\Omega_t) = -342144 t^3 + O(t^4).
\end{equation}
This quantity is negative for all sufficiently small
$t>0$.  Geometrically, any line of negative slope through
the origin, sufficiently close to the horizontal,
contains points of $(f,g)$ on both sides.

Another edge of $X_{\beta}$ has endpoints
$(8,0,0,8,8,0)$ and $(4,4,4,4,4,4)$.
We restrict our functions to the point
\begin{equation}
\Psi_t=(1-t^2) (8,0,0,8,8,0) + t^2 (4,4,4,4,4,4).
\end{equation}
a point which lies in an edge of $X_{\beta}$ for
small $t$.
Using Mathematica, we compute that  
Let \begin{equation}
(12-t)g(\Psi_t)-f(\Psi_t)=-57344 t^5 + O(t^6).
\end{equation}
This function is negative for all sufficiently
small $t>0$.
Geometrically, if we take the line of slope
$1/12$ through the origin and increase
its slope by any small positive amount, 
points of $(f,g)$ will lie on either side of
the line.

\subsection{A Pair of Incident Edges}
\label{incidentpair}

Let $\beta$ be a pair of incident edges.  Let
$X_{\beta}$ is the set of $4$ chambers $X_D$ such
that $\beta$ is disjoint from the outer two
edges of $D$ and the black dot is the common
endpoint of the two edges of $\beta$.
Put another way, the shortest vertex-sum occurs
at the vertex incident to the two edges of $\beta$,
and the pair of opposites with the largest axis sum
is disjoint from $\beta$ except at the endpoints.
There are $4$ such chambers.  Figure 5.2 shows
the decorations corresponding to $2$ of the
chambers. These decorations correspond to
the simplices listed below.
The other $2$ are the images of
these two under the element of $K_4$ which
stabilizes $\beta$.

\begin{center}
\resizebox{!}{1.5in}{\includegraphics{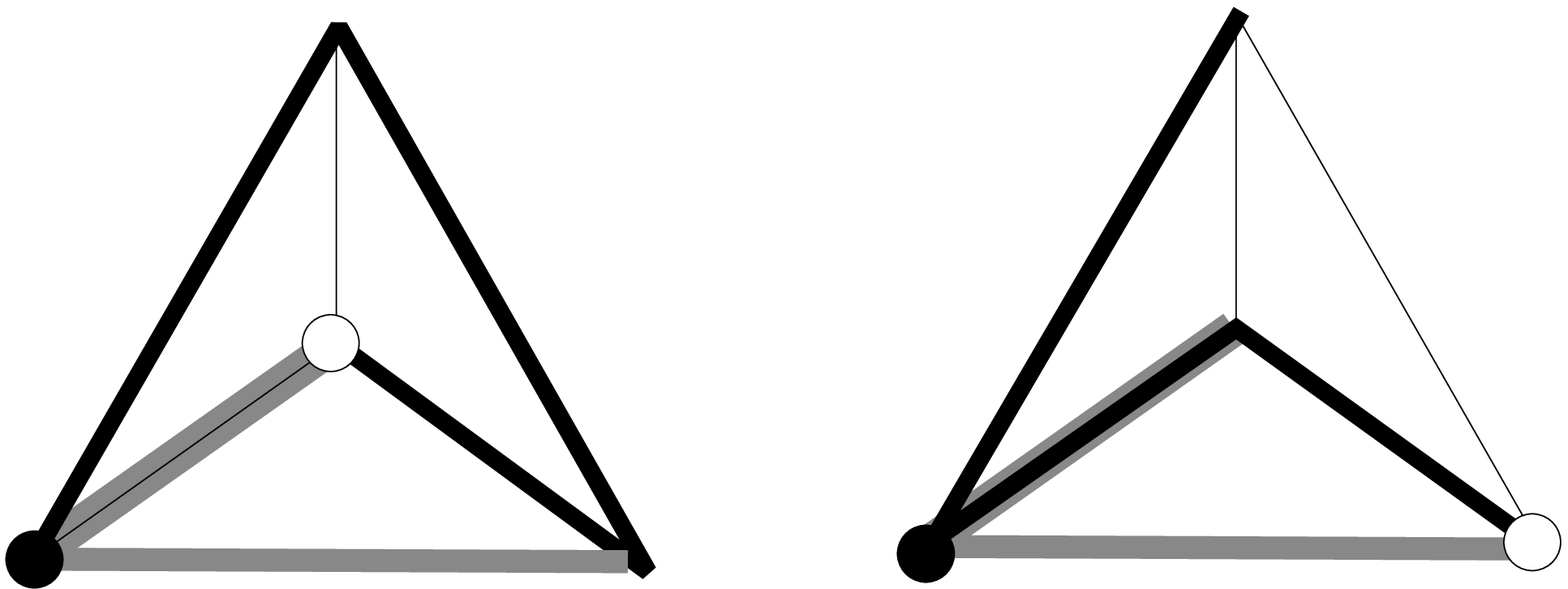}}
\newline
{\bf Figure 5.2:\/} Two of four the decorations for $X_{\beta}$.
\end{center}

We take $\beta=\{e_{12},e_{13}\}$. First, we
anti-certify all the chambers not in $X_{\beta}$.
Now we turn to $X_{\beta}$.

Two of the chambers
of $X_{\beta}$ intersect $X_{24}$ in the simplices
$$(4,4,4,4,4,4),\ (0,8,0,8,0,8),\ (4,0,4,4,8,4)$$
$$(3,3,6,6,3,3),\ (8,0,0,8,8,0),\ (0,6,6,6,6,0)$$
and
$$(4,4,4,4,4,4),\ (0,8,0,8,0,8),\ (4,0,4,4,8,4)$$
$$(3,3,6,6,3,3),\ (0,8,0,8,0,8),\ (0,6,6,6,6,0)$$
The other two simplices are images of these
under the order $2$ subgroup of $S_4$ which
stabilizes $\beta$.

We work with the functions
\begin{equation}
P=g \hskip 30 pt Q=2g - f.
\end{equation}
The PDA takes $421$ steps to certify that $P \geq 0$
on each simplex, and $479$ steps to certify that
$Q \geq 0$ on each simplex.
These calculations show that $A_{\\beta} \leq 0$
and $B_{\beta} \geq 12$.  Now we prove equality.

Let $V_1,....,V_6$ be the $6$ vertices of the
first simplex listed above. Call this simplex $\Sigma$. Define
\begin{equation}
\Omega_t=(1-t-t^2)V_5+ t V_2 + t^2 V_3 \in \Sigma.
\end{equation}
We compute that
\begin{equation}
g(\Omega_t)+t(\Omega_t)=-2097152 t^7 + O(t^8).
\end{equation}
We also compute that
\begin{equation}
(2-t)vf(V_1) - f(V_1)=-8192t.
\end{equation}
The same argument as in the previous section
shows that $(A_{\beta},B_{\beta})=(0,12)$.

\subsection{A Pair of Opposite Edges}

Let $\beta$ be a pair of opposite edges.  Let
$X_{\beta}$ be the set of $32$ chambers
$X_D$ such that $\beta \not \subset D$.
In other words, the largest axis sum
does not occur at $\beta$.

We take $\beta=\{e_{12},e_{34}$.
First, we anti-certify the chambers of
$X-X_{\beta}$.  

Using the symmetry of the permutation group,
it suffices to consider the $4$ chambers which
intersect $X_{24}$ in the following simplices.

$$(4,4,4,4,4,4),\ (0,8,0,8,0,8),\ (4,0,4,4,8,4)$$
$$(3,6,3,3,6,3),\ (8,0,0,8,8,0),\ (0,6,6,6,6,0)$$
and
$$(4,4,4,4,4,4),\ (0,8,0,8,0,8),\ (4,0,4,4,8,4)$$
$$(3,6,3,3,6,3),\ (8,0,0,8,8,0),\ (6,6,0,0,6,6)$$
and
$$(4,4,4,4,4,4),\ (8,8,8,0,0,0),\ (4,0,4,4,8,4)$$
$$(3,6,3,3,6,3),\ (8,0,0,8,8,0),\ (0,6,6,6,6,0)$$
and
$$(4,4,4,4,4,4),\ (8,8,8,0,0,0),\ (4,0,4,4,8,4)$$
$$(3,6,3,3,6,3),\ (8,0,0,8,8,0),\ (6,6,0,0,6,6)$$

Figure 5.3 shows the corresponding decorations.

\begin{center}
\resizebox{!}{4.5in}{\includegraphics{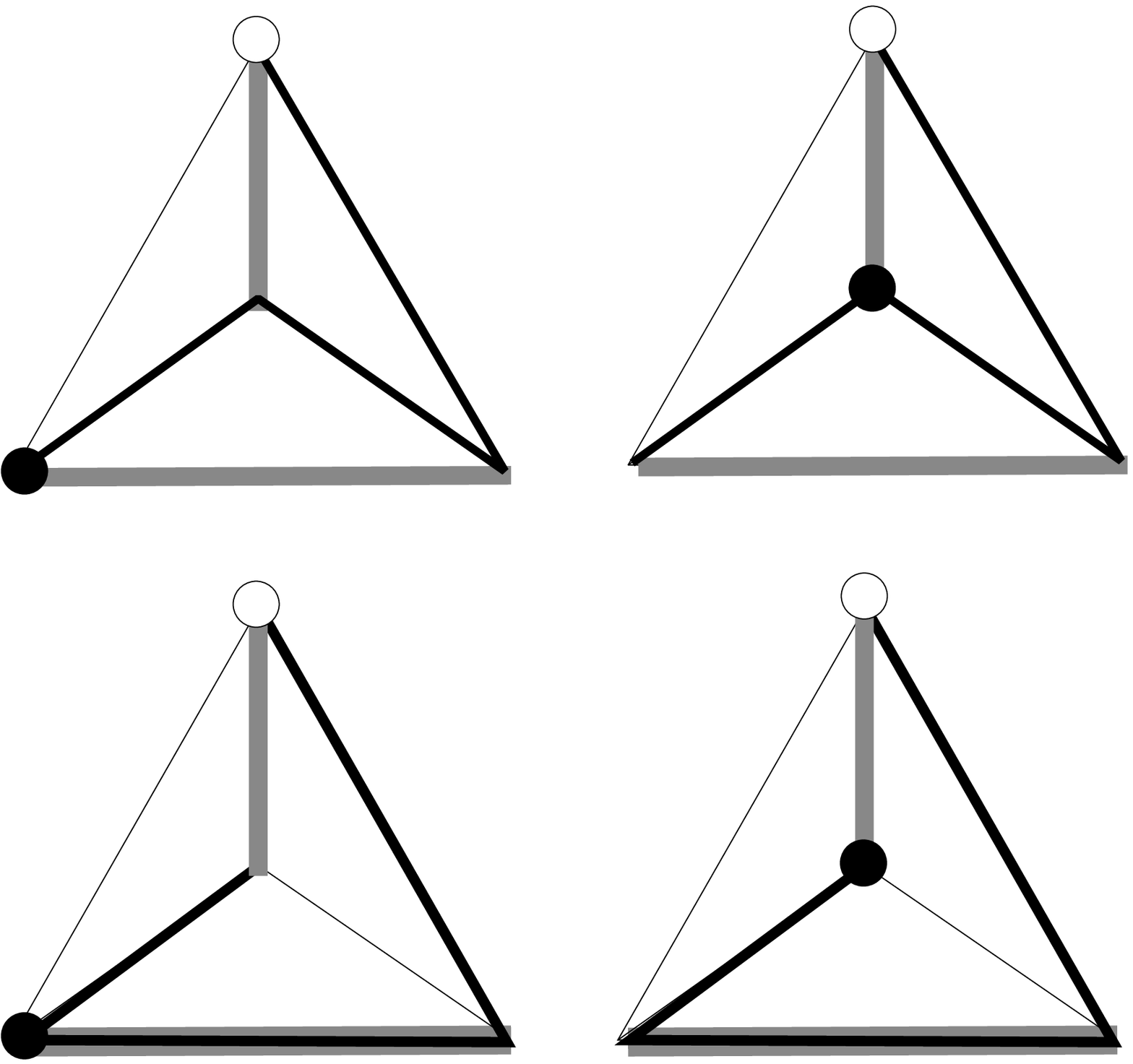}}
\newline
{\bf Figure 5.3:\/} Three of the $12$ decorations for $X_{\beta}$.
\end{center}

We work with the functions
\begin{equation}
P=g, \hskip 30 pt Q=6g-f,
\end{equation}

The Positive Dominance algorithm shows that $P \geq 0$ on
the above simplices in $473,473,331,331$ steps respectively.
The Positive Dominance algorithm shows that $Q \geq 0$ on
the above simplices in $467,467,1161,1161$ steps respectively.
Our calculations show that $A_{\beta} \leq 0$ and
$B_{\beta}\geq 4$. 
\newline
\newline
In fact, $A_{\beta}=0$ and $B_{\beta} \in (4,5)$.
I don't know the precise value of $B_{\beta}$ and
I'm not sure that the PDA could establish it
even if I knew what it was.

\subsection{Tripods}

Let $\beta$ be a triple of edges incident to the
same vertex $v$.  Let
$X_{\beta}$ be the set of $12$ chambers
$X_D$ such that the black vertex is $v$.
That is, the apex of the tripod has the
smallest vertex-sum.

We take $\beta=\{e_{12},e_{13},e_{14}\}$. 
First of all, we anti-certify all the
chambers of $X-X_{\beta}$.  Now we turn
to $X_{\beta}$.

By symmetry, it suffices to consider 
the first and third simplex listed in the
previous section.
Both these simplices are contained in
the simplex $\C_{21}$ from the $12$-partition
discussed in \S 2.  We will work with
$\C_{21}$ because it is just a single simplex.
The vertices of $\C_{21}$ are

$$(4,4,4,4,4,4),\ (0,8,0,8,0,8),\ (8,0,0,8,8,0)$$
$$(0,0,8,0,8,8),\ (0,6,6,6,6,0),\ (6,6,0,0,6,6)$$

We work with the functions
\begin{equation}
P=4g-3f, \hskip 30 pt
Q=3g-f.
\end{equation}
The PDA shows that $P \geq 0$ on $\C_{21}$ in $967$ steps.
The PDA shows that $Q \geq 0$ on $\C_{21}$ in $779$ steps.
These calculations show that $A_{\beta} \leq 8$ and
$B_{\beta} \geq 18$.  

We compute
\begin{itemize}
\item $f(8,8,8,8,8,8)=3 \times 2^{12}$.
\item $g(8,8,8,8,8,8)=2^{12}$.
\item $f(6,6,0,0,6,6)=-2^7 \times 3^6$.
\item $g(6,6,0,0,6,6)=-2^7 \times 3^5$.
\end{itemize}
These calculations show that the two images
of $(f,g)$ lie on the two boundary rays of
the cones defined by the conditions $P \geq 0$
and $Q \geq 0$.   Hence
$(A_{\beta},B_{\beta})=(8,18)$.

\subsection{$3$-paths}

Let $\beta$ be a $3$-path.  Let
$X_{\beta}$ be the set of $8$ chambers
$X_D$ such that the black vertex is an
interior vertex of $\beta$, and
the outer edges of $D$ are disjoint from
$\beta$. In other words, the smallest
vertex-sum occurs at an interior vertex
of $\beta$ and the largest axis-sum occurs
on a pair of opposite edges disjoint from
$\beta$ (except at the endpoints.)

We take $\beta=\{e_{12},e_{14},e_{23}\}$.
First of all, we anti-certify all the
chambers of $X-X_{\beta}$.  Now we turn
to $X_{\beta}$.
The chambers of $X_{\beta}$ all lie in
$\C_{21} \cup \C_{13}$.  By symmetry, it suffices
to consider the chambers in $\C_{21}$, the simplex
from the previous section.
we work with the functions
\begin{equation}
P=4g+f, \hskip 30 pt Q=3g-3f.
\end{equation}
The PDA certifies that $P \geq 0$ in $\C_{21}$ in
$823$ steps.  The PDA certifies that
$Q \geq 0$ on $\C_{21}$ in $1243$ steps.
These calculations show that
$A_{\beta} \leq -6$ and $B_{\beta} \geq 16$.
\newline
\newline
{\bf Remark:\/}
I don't know the optimal constants in this case.

\subsection{$4$-Cycles}

Let $\beta$ be a $4$-cycle.  Let $X_{\beta}$ be
the $16$ chambers $X_D$ such that the outer two edges
of $D$ are disjoint from $\beta$.  In other words,
the axis-sum of the $\beta$-complement is largest.

We take 
$\beta=\{e_{12},e_{13},e_{24},e_{34}\}$.
First of all, we anti-certify the chambers of $X-X_{\beta}$.
Turning to $X_{\beta}$, we work with the polynomials
\begin{equation}
P=g, \hskip 30 pt Q=g-f.
\end{equation}

In this case, all the chambers of $X_{\beta}$
intersect $X_{24}$ inside the simplex
$\A_3$.
By symmetry, it suffices to prove that
$P \geq 0$ and $Q \geq 0$ on the simplex
$\C_{31}$.  This simplex has vertices

$$(4,4,4,4,4,4),\ (0,0,8,0,8,8),\ (8,0,0,8,8,0)$$
$$(0,8,0,8,0,8),\ (0,6,6,6,6,0),\ (6,0,6,6,0,6)$$

The PDA takes $755$ steps to verify that
$P \geq 0$ on $\C_{31}$ and
$1687$ steps to verify that
$Q \geq 0$ on $\C_{31}$.
These calculations show that $A_{\beta} \leq 0$ and
$B_{\beta} \geq 24$.

Now we observe that
$$f(4,4,4,4,4,4)=g(4,4,4,4,4,4)=2^{14}.$$
This forces $B_{\beta}=24$.

Let
$V_1,...,V_6$ be the vectors listed above and
we define
\begin{equation}
\Theta_t=(1-t-t^2)V_2 + t V_3 + t^2 V_4.
\end{equation}
We compute that
\begin{equation}
(g+tf)(\Theta_t)=-8388608 t^7 + O(t^8).
\end{equation}
This shows that $A_{\beta}=0$.

\subsection{$3$-cycles}

We have now completes the proof of
Theorem \ref{schema} and all the auxiliary
facts mentioned after that result in the
introduction.  Now we turn to the proof
of Theorem \ref{schema3}.

Let $\beta$ be a $3$-cycle in $K_4$. Let
$X_{\beta}$ be the $36$ chambers
$X_D$ so that the black vertex lies in $\beta$.
That is, the vertex with the shortest vertex-sum
must be a vertex of $\beta$.   We take
$\beta=\{e_{12},e_{13},e_{23}\}$.  We first
anti-certify the chambers of $X-X_{\beta}$.

Now we turn to $X_{\beta}$.
By symmetry,
it suffices to consider points of $X_{24}$
where the vertex-sum is smallest at vertex $1$.
But then we are talking about points in
the simplex $\B_1$.  This simplex has vertices

$$(0,6,6,6,6,0),\ (6,0,6,6,0,6),\ (6,6,0,0,6,6)$$
$$(8,0,0,8,8,0),\ (0,8,0,8,0,8),\ (0,0,8,0,8,8)$$

We work with the function
\begin{equation}
P=3g-f.
\end{equation}
The PDA certifies that $P \geq 0$ on $\B_1$ in $1275$ steps.

Let $V_1,...,V_6$ be the vertices listed above.
Define
\begin{equation}
\Omega_t=(1-t^2)V_1 + t^2 V_2, \hskip 30 pt
\Psi_t=(1-t-t^2)V_4 + t V_1 + t^2 V_2.
\end{equation}
We compute that
\begin{equation}
(3+t)g(\Omega_t)-f(\Omega_t)=-497664 t^5 + O(t^6).
\end{equation}
Geometrically, this means that when we rotate the line
$3y=x$ about the origin
in such a way as to  slightly increase its slope,
the upper half plane bounded by the new line does
not contain the image $\B_1$ under $(f,g)$.
We compute that
\begin{equation}
(3-t)g(\Psi_t)-f(\Omega_t)=663552 t^6 + O(t^7).
\end{equation}
Geometrically, this means that when we rotate the line
$3y=x$ about the origin in such a way as to slightly decrease its slope,
the upper half plane bounded by the new line does
not contain the image $\B_1$ under $(f,g)$.
This shows that $g-Cf \geq 0$ on $X_{\beta}$ if and only if
$C=3$.  The corresponding statement in Theorem \ref{schema3}
follows from this fact, and from homogeneity.

This completes the proof of
Theorem \ref{schema3}.

\subsection{The Unfriendly Configurations}

Aside from a $3$-cycle, there are two unfriendly
configurations, $\alpha$ and $\beta$, chosen
so that $$K_4-\alpha=\{e_{12},e_{13}\}, \hskip 30 pt
K_4-\beta=\{e_{12}\}.$$

Some experimental evidence suggests that Corollary
\ref{schema2} holds for the set $X_{\alpha}$
of $24$ chambers of the form $X_D$, where
$\alpha \not \subset D$ and the black dot is
incident to at least $2$ edges of $\alpha$. 
However, computer plots also show that the
image of any chamber $X_D$ under the
map $(f,g)$ does not lie in a halfspace.
(Nonetheless it seems that $f>0$ implies $g>0$.)
So, our method simply does not apply here.

Some experimental evidence suggests that
Corollary \ref{schema2} holds for a certain
set $X_{\beta}$ of $24$ chambers whose
description in terms of the decorations is
rather complicated.  $16$ of the chambers
intersect $X_{24}$ inside $\A_1$, and the
remaining $8$ intersect $X_{24}$ inside
$(\B_3 \cup \B_4)-\A_1$. The interested
reader can download our program and see the
set exactly.  

Sitting inside $X_{\beta}$ is
a smaller set $X'_{\beta}$ consisting
of the $8$ chambers which intersect
$X_{24}$ inside $\C_{13} \cup \C_{14}$.
When $X_D$ is one of these $8$ chambers,
the image of $X_D$ under
$(f,g)$ is contained in a proper cone in
$\R^2$, and the volume-increase part of
Corollary \ref{schema2} holds.  We also
leave this to the interested reader.

\newpage

%% file: 6exist.tex
\section*{Appendix: Existence in all Dimensions}

Here I'll give the proof of
Theorem \ref{exist}. I learned all the
arguments here from Peter Doyle and Igor Rivin.
Let $D$ stand for
a list $\{d_{ij}\}$.  We will perform
operations componentwise, so that
$D+t=\{d_{ij}+t\}$, etc.
We define tetrahedral lists in all dimensions just
as in the $3$ dimensional case. Let
$T$ denote the space of tetrahedral lists.

\begin{theorem}
\label{sq}
If $D \in T$ then
$\sqrt D \in T$.
\end{theorem}

See [{\bf WW\/}, Corollary 4.8].  In [{\bf WW\/}] this
result is attributed to Von Neumann, though Rivin
calls it Schoenberg's result.

\begin{theorem}
\label{conv}
If $A,B \in T$ then
$\sqrt{A^2+B^2} \in T$.
The simplex represented by
$\sqrt{A^2+B^2}$ has larger volume than the
simplex represented by $A$.
\end{theorem}

\startproof
See [{\bf R\/}] for a proof. 
See also [{\bf BC\/}, Lemma 1].  Here is a
self-contained proof.

Given a quadratic form $Q$ and
a linear isomorphism $L$, we have the new quadratic form
\begin{equation}
L^*Q(v,w)=Q(L^{-1}(v),L^{-1}(w)).
\end{equation}
 Let $\Delta \subset \R^n$ denote some copy
of the regular simplex.
Let $L_A$ denote the linear transformation which carries
$\Delta$ to a simplex $\Delta_A$ whose lengths are realized by $A$.
Let $Q_A$ be the quadratic form such that
$L_A^*(Q_A)$ is the standard quadratic form -- i.e. the
dot product.  By construction, $Q_A$ assigns the length
list $A$ to the sides of $\Delta$.  Likewise define $Q_B$.
The sum $Q_C=Q_A+Q_B$ is also positive definite.
Given an edge $e$ of $\Delta$.   Let $C$ denote the
list of lengths that $Q_C$ assigns to $\Delta$. We compute
$$C_e=\sqrt{Q_C(e,e)}=\sqrt{Q_A(e,e)+Q_B(e,e)}=\sqrt{A_e^2+B_e^2}.$$
This shows that $Q_C$ assigns the corresponding number on the
list  $\sqrt{A^2+B^2}$ to the edge $e$.   But then let
$L_C$ be the linear transformation which pushes $L_C$
forward to the standard quadratic form.  By construction
$L_C(\Delta)$ is the simplex realizing the list
$\sqrt{A^2+B^2}$.

Note that $L_C(e)>L_A(e)$ for all $e \in \R^n$.  But then
the linear map carrying $\Delta_A$ to $\Delta_C$ strictly
increases all distances, and hence also increases volume.
\endproof

Let $\phi_t$ denote the lengthening flow. That is,
$\phi_t(D)=D+t$.  Let $D_t=\phi_t(D)$.

\begin{lemma}
\label{tangent}
$$\frac{d}{dt}{\rm vol\/}(D_t) \geq 0$$
for all $D \in T$.
\end{lemma}

\startproof
Choose some point $D \in T$.
By Theorem \ref{sq} and scaling,
$\sqrt{2tD} \in T$.  By Theorem \ref{conv}.
$\sqrt{D^2+2tD} \in T$.
But we can write out
\begin{equation}
\sqrt{D^2+2tD}=D(\sqrt{1+2t/D}) = D+t + {\rm higher\ order\ terms\/}
\end{equation}
Theorem \ref{conv} shows that
$${\rm vol\/}(D+t+{\rm higher\ order\ terms\/})>{\rm vol\/}(D_0).$$
Taking the limit as $t \to 0$ we get
the result of this lemma.
\endproof

Lemma \ref{tangent} implies that $D_t \in T$ for all $t$,
because ${\rm vol\/}(D_t)$ cannot converge to $0$.
This takes part of the existence statement in
Theorem \ref{exist}.

To see that ${\rm vol\/}(D_t)$ is strictly increasing, we
look more carefully at the proof of Theorem \ref{conv}.
When all the terms on the list $B$ have size $O(t)$, and
$t$ is much smaller than the minimum length on the list $A$,
the quadratic form $Q_C$ assigns the length
$A_e+O(t)$ to each unit length vector $e \in \R^n$.
But then the linear map $\Delta_A$ to $\Delta_C$ increases all
unit distances by $O(t)$, and hence increases volume by $O(t)$ as well.
This gives
$${\rm vol\/}(D+t+{\rm higher\ order\ terms\/})>{\rm vol\/}(D_0)+C_Dt,$$
for some positive constant $C_D$.
Taking the limit, we see that ${\rm vol\/}(D_t)$ is strictly
increasing, and in fact has positive derivative.
This completes the proof of Theorem \ref{exist}.

\newpage

%% file: refs.tex
\section{References}

[{\bf BC\/}] K Bezdek and R. Connelly, 
{\it Pushing Disks apart -- The Kneser-poulsen Conjecture in the plane\/},
Journal fur die riene und angewandte Mathematik (2002) pp 221-236
\newline
\newline
[{\bf L\/}] F. Luo, {\it $3$-Dimensional Schlafli Formula and its Generalizations\/}, \newline
arXiv:0802.2580
\newline
\newline
[{\bf P\/}] I. Pak, {\it The Cayley-Menger Determinant\/},
\newline
UCLA course notes (2006) \newline
{\bf http://www.math.ucla.edu/$\sim$pal/courses/geo/cm.pdf\/}
\newline
\newline
[{\bf Sa\/}] I. Kh. Sabitov, {\it Algebraic methods for solution of polyhedra\/}, Russian Math Surveys {\bf 66.3\/} (2011) pp 445-505
\newline
\newline
[{\bf R\/}] I. Rivin, {\it Some Observations on the Simplex\/}, preprint (2003) arXiv 0308239
\newline
\newline
[{\bf S\/}] R. Schwartz,
{\it The Projective Heat Map on Pentagons\/}, \newline
research monograph (2014) preprint
\newline
\newline
[{\bf W\/}], S. Wolfram, {\it Mathematica\/}, Wolfram Media and Cambridge University Press (1999)
\newline
\newline
[{\bf WW\/}] J. H. Wells and L. R. Williams, {\it Embeddings and Extensions in Analysis\/},
Springer-Verlag 1975